\documentclass[11pt,reqno,a4paper]{amsart}

\usepackage{etoolbox}
\usepackage{amsmath}
\usepackage{enumerate}
\usepackage{amssymb}
\usepackage{amscd}
\usepackage{amsthm}
\usepackage{amsfonts}
\usepackage{graphicx}
\usepackage[all,cmtip]{xy}
\usepackage{enumitem}

\patchcmd{\subsection}{-.5em}{.5em}{}{}
\patchcmd{\subsubsection}{-.5em}{.5em}{}{}

\usepackage{enumitem}

\usepackage[T1]{fontenc}

\usepackage[scaled]{helvet} 
\usepackage{courier} 
\usepackage{eulervm}
\normalfont

\makeatother
\usepackage{hyperref}

\numberwithin{equation}{section}

\newcommand{\SL}{\operatorname{SL}}

\newcommand{\cL}{\mathcal{L}}

\newcommand{\cP}{\mathcal{P}}
\newcommand{\cQ}{\mathcal{Q}}

\newcommand{\cY}{\mathcal{Y}}

\newcommand{\bR}{\mathbb{R}}

\newcommand{\qen}{\enskip \textrm{and} \enskip}
\newcommand{\qand}{\quad \textrm{and} \quad}

\newcommand\subsetsim{\mathrel{%
\ooalign{\raise0.2ex\hbox{$\subset$}\cr\hidewidth\raise-0.8ex\hbox{\scalebox{0.9}{$\sim$}}\hidewidth\cr}}}
\newcommand{\eps}{\varepsilon}

\DeclareMathOperator{\cum}{Cum}

\theoremstyle{theorem}
\newtheorem{theorem}{Theorem}[section]
\newtheorem{corollary}[theorem]{Corollary}
\newtheorem{proposition}[theorem]{Proposition}
\newtheorem{lemma}[theorem]{Lemma}

\theoremstyle{definition}

\begin{document}

\title{Central Limit Theorems in the Geometry of Numbers}

\author{Michael Bj\"orklund}
\address{Department of Mathematics, Chalmers, Gothenburg, Sweden}
\email{micbjo@chalmers.se}

\author{Alexander Gorodnik}
\address{University of Bristol, Bristol, UK}
\email{a.gorodnik@bristol.ac.uk}



\date{Central limit theorems, Diophantine approximation}


\maketitle

\begin{abstract}
We investigate in this paper the distribution of the discrepancy of various lattice counting functions.
In particular, we prove that the number of lattice points contained in certain domains
defined by products of linear forms satisfies a Central Limit Theorem.
Furthermore, we show that the Central Limit Theorem holds 
for the number of rational approximants for weighted Diophantine approximation in $\mathbb{R}^d$.
Our arguments exploit chaotic properties of the Cartan flow on the space of lattices.
	
\end{abstract}

\section{Introduction}

Let $\Omega_T\subset \mathbb{R}^d$ be an increasing family of compact domains,
and let $\cL_d$ denote the space of lattices in $\bR^d$ with covolume one, endowed with the unique 
$SL_d(\bR)$-invariant probability measure $\lambda_d$. We consider the counting function $\Lambda \mapsto 
|\Lambda\cap \Omega_T|$ on $\cL_d$. Under mild assumptions on the domains $\Omega_T$, 
$$
|\Lambda\cap \Omega_T|\sim \frac{\hbox{vol}(\Omega_T)}{\hbox{vol}(\mathbb{R}^d/\Lambda)}\quad \textrm{as $T\to\infty$, for all $\Lambda \in \cL_d$.}
$$
One may ask whether it is possible to derive more precise information about the asymptotic behavior of $|\Lambda\cap \Omega_T|$ for \emph{generic} lattices $\Lambda$. \\

The following naive heuristics might give an idea of what to expect. Let us decompose 
$$
\Omega_T=\bigsqcup_{i=1}^N \Omega_T^{(i)}
$$
into regions $\Omega_T^{(i)}$ with $\hbox{vol}(\Omega_T^{(i)})\approx 1$. Then
$$
|\Lambda\cap \Omega_T|=\sum_{i=1}^N |\Lambda\cap \Omega^{(i)}_T|,
$$
and provided that $\Omega^{(i_1)}$ and $\Omega^{(i_2)}_T$ are "far apart", it seems plausible to conjecture that
the random variables $\Lambda \mapsto |\Lambda \cap \Omega^{i_j}_T|$, for $j = 1,2$, on $\cL_d$ are "almost independent".
Thus, one might wonder whether $|\Lambda\cap \Omega_T|$ behaves like a sum of independent random variables. 

Some classical results of W. Schmidt motivated our line of study. In \cite{sch1, sch2}, Schmidt showed that for generic lattices $\Lambda \in \cL_d$,
$$
|\Lambda\cap \Omega_T|= \frac{\hbox{vol}(\Omega_T)}{\hbox{vol}(\mathbb{R}^d/\Lambda)}
+O_{\Lambda,\eps}\left(\hbox{vol}(\Omega_T)^{\frac{1}{2}+\eps}\right)\quad \hbox{ for all $\eps>0$.},
$$
and thus the counting function $\Lambda \mapsto |\Lambda\cap \Omega_T|$ indeed exhibits 
cancellations of the same order as a sum of independent random variables. Remarkably, the argument 
in \cite{sch1} implicitly follows the heuristic approach outlined above and proves some form of pairwise 
independence using arithmetic considerations.

The aim of this work is to establish a Central Limit Theorem (CLT) in this setting, at least under some additional assumptions on 
the domains $\Omega_T$. We stress that it is unlikely that a CLT holds for general domains;  for instance, the counting of lattice 
points in the regions 
$$
\Omega_T = \big\{ (x,y) \in \mathbb{R}^2 \, : \,  0 < x < Ty \qand 1 < y < 2 \big\}
$$
is closely related to the distribution of averages for the horocyclic flow on $\cL_2$, which do not admit a CLT (see e.g. \cite{FF}).

We shall in this paper consider domains defined by products of linear forms on $\bR^d$. Such domains can be tesselated 
using images of a small number of regular tiles under a family of diagonal matrices in $\SL_d(\bR)$, and allows us to use 
dynamical  arguments developed in our recent work \cite{BG}.  The crucial ingredients in our approach are quantitative estimates
on higher-order correlations established in our joint work with Einsiedler \cite{BEG}. Besides generic lattices in the space of 
lattices $\cL_d$, we also consider the family of lattices
\begin{equation}
\label{eq:lll}
\Lambda_{\overline{x}}=\{(p_1-q x_1,\ldots, p_d-qx_d, q):\, (\overline{p},q)\in \mathbb{Z}^d\times \mathbb{Z}\},
\end{equation}
for $\overline{x}=(x_1,\ldots,x_d)\in\mathbb{R}^d$, which arises in many problems in the theory of Diophantine approximation.

\section{Main results}

\subsection{Distribution of values for products of linear forms}

We fix a collection of linearly independent linear forms
$L_1,\ldots,L_d:\mathbb{R}^d\to\mathbb{R}$ with $d\ge 3$
and consider the product form 
$$
N(x)=L_1(x)\cdots L_d(x).
$$
Our aim is to analyze the distribution of the values $N(x)$ when $x$ belongs to a lattice in $\mathbb{R}^d$.
We fix an interval $(a,b)\subset \mathbb{R}^+$, and for $T\ge 1$, we define the domains
\begin{equation}
\label{eq:omega1}
\Omega_T=\big\{x \in \mathbb{R}^d \, : \, N(x)\in (a,b) \qand |L_1(x)|,\ldots, |L_d(x)|< T\big\}.
\end{equation}
It is not hard to show that 
$$
\hbox{vol}(\Omega_T)=c\,(b-a)(\log T)^{d-1}+O\left((\log T)^{d-2}\right),
$$
for some $c=c(L_1,\ldots,L_d)>0$, and by \cite{sch2}, almost all unimodular lattices $\Lambda$ in $\mathbb{R}^d$ satisfy
$$
|\Lambda\cap \Omega_T| = \hbox{vol}(\Omega_T)+O_{\Lambda,\eps}\left((\log T)^{\frac{d-1}{2}+\eps}\right)\quad \hbox{for all $\eps>0$.}
$$
We shall investigate how the error term (also known as the discrepancy) in this formula behaves. \\

Our first result shows that the error term admits a Central Limit Theorem. 
We have currently verified our argument for $d\ge 4$, but it might be possible optimize the estimates to deal with the case $d=3$ as well.

\begin{theorem}[CLT for lattice counting]\label{th:lattice1}
For $d\ge 4$, there exists explicit $\sigma>0$ such that
for every $u\in \mathbb{R}$,
$$
\lambda_d\left(\left\{\Lambda\in \mathcal{L}_d:\, \frac{|\Lambda\cap \Omega_T| - \hbox{\rm vol}(\Omega_T)}{\hbox{\rm vol}(\Omega_T)^{1/2}}<u \right\} \right)
\longrightarrow
\frac{1}{\sqrt{2\pi\sigma}} \int_{-\infty}^u e^{-t^2/(2\sigma)}\, dt
$$
as $T\to\infty$.
\end{theorem}

A more general version of this theorem can be established along similar lines. Instead of considering linear forms, let $L_i:\mathbb{R}^d\to \mathbb{R}^{d_i}$, $1\le i\le k$, be a family of linear maps, and set
$$
N(x)=\prod_{i=1}^k \|L_i(x)\|^{d_i},
$$
where $\|\cdot\|$ denotes the Euclidean norm. We shall assume that $d_1+\cdots +d_k=d$ and the map $(L_1,\ldots,L_k)$ defines a
bijection $\mathbb{R}^d\to \prod_{i=1}^k \mathbb{R}^{d_i}$. For a fixed interval $(a,b)\subset \mathbb{R}^+$, we define the domains
\begin{equation}
\label{eq:omega2}
\Omega_T=\big\{ x \in \mathbb{R}^d \, : \,  N(x)\in (a,b) \qand  \|L_1(x)\|,\ldots, \|L_k(x)\|< T\}.
\end{equation}
We show that a version of Theorem \ref{th:spiralling} still holds for these domains. In the case $k = 2$, such a result was 
also established in \cite{dfv}, using a different method, which does not seem to generalize to $k \geq 3$.

\subsection{Spiralling}

Motivated by the paper \cite{agt}, we shall also study  ``spiraling'' of the lattice points
contained in the regions \eqref{eq:omega2}, that is, the distribution 
of their angular components.
We denote by $\omega_i:\mathbb{R}^{d_i}\backslash \{0\}\to S^{d_i-1}$ the radial projections,
and for a lattice $\Lambda$ in $\mathbb{R}^d$ and a Borel set $D \subset S$, we define
\[
S_T(\Lambda,D) = \big\{ x \in \Omega_T \cap \Lambda \, : \, (\omega_1(x),\ldots,\omega_k(x)) \in D \big\}.
\]
It is not hard to show (see \cite{agt}) that for almost every unimodular lattice $\Lambda$ in $\mathbb{R}^d$,
and for any Borel subset $D\subset S$,
$$
\frac{|S_T(\Lambda,D)|}{\hbox{\rm vol}(\Omega_T)} \to \hbox{vol}(D)\quad\hbox{ as $T\to\infty$}.
$$
We prove that if some regularity is imposed on $D$, then a suitable Central Limit Theorem also holds. 
Currently, the argument has been verified for $d\ge 4$, but it might be possible optimize the estimates further 
to deal with the case $d=3$.

\begin{theorem}[CLT for spiraling]
	\label{th:spiralling}
	For $d\ge 4$ and for every domain $D \subset S$ with piecewise smooth boundary, there exists explicit $\sigma=\sigma(D)>0$ such that for every $u\in \mathbb{R}$,
	$$
	\lambda_d\left(\left\{\Lambda\in \mathcal{L}_d:\, \frac{|S_T(\Lambda,D)| - \hbox{\rm vol}(D) \hbox{\rm vol}(\Omega_T)}{\hbox{\rm vol}(\Omega_T)^{1/2}}<u \right\} \right)
	\longrightarrow
	\frac{1}{\sqrt{2\pi\sigma}} \int_{-\infty}^u e^{-t^2/(2\sigma)}\, dt
	$$
	as $T\to\infty$.
\end{theorem}

\subsection{Diophantine approximation}

Let us now discuss the distribution of integral solutions of some inequalities which arise in the theory
of Diophantine approximation. We start with Diophantine approximation on the real line, which is better 
understood due to the theory of continued fractions. Fix $c > 0$, and for $x\in \mathbb{R}$, we consider 
the Diophantine inequality
\begin{equation}
\label{eq:dioph}
\left|x-\frac{p}{q}\right|< \frac{c}{q^2}
\end{equation}
with $(p,q)\in \mathbb{Z}\times \mathbb{N}$, and the corresponding counting function
$$
N_T(x)=|\big\{(p,q)\in \mathbb{Z}\times \mathbb{N}:\, 1\le q\le T \qen \hbox{$\frac{p}{q}$ is a solution of \eqref{eq:dioph}}  \big\}|
$$
It is known (see, for instance, \cite{sch1}) that for almost every $x\in [0,1]$,
$$
N_T(x)= 2c\, \log T+O_{x,\eps}((\log T)^{1/2+\eps}), \quad \hbox{for all $\eps>0$.}
$$
Fuchs showed in \cite{f} that the discrepancy in this formula satisfies the Central Limit Theorem, that is to say, 
there exists $\sigma>0$ such that for every $u\in \mathbb{R}$,
\begin{equation}
\label{eq:CLT1}
\left|\left\{x\in [0,1]:\, \frac{N_T(x)-2c\,\log T}{(\log T\cdot\log\log T)^{1/2}}<u\right\}\right|
\longrightarrow \frac{1}{\sqrt{2\pi\sigma}} \int_{-\infty}^u e^{-t^2/(2\sigma)}\,dt
\end{equation}
as $T\to\infty$. We stress that the correct normalization in \eqref{eq:CLT1} has caused some confusion in the previous works
\cite{lev1,lev2,ph}; the additional $(\log\log T)$-factor arises here because a certain counting function on $\cL_2$ is \emph{not} 
square-integrable. This non-integrability issue does not appear in higher dimensions, whence this additional normalization 
factor should disappear. An analogue of this result for simultaneous Diophantine approximation has been recently established in \cite{dfv}. 

In this paper we consider the following more general problem in \emph{weighted} Diophantine approximation.
Let us fix a collection of weights 
$$
0<w_1,\ldots, w_d<1 \qand  w_1+\ldots+w_d=1,
$$
and constants $c_1,\ldots,c_d > 0$. Given a vector $\overline{x}=(x_1,\ldots,x_d)\in \mathbb{R}^d$,
we are interested in understanding the asymptotics of solutions for the system of Diophantine inequalities
defined by
\begin{equation}
\label{eq:diop2}
\left|x_1-\frac{p_1}{q}\right|< \frac{c_1}{q^{1+w_1}},\;\;\ldots,\;\;
\left|x_d-\frac{p_d}{q}\right|< \frac{c_d}{q^{1+w_d}}
\end{equation}
with $(\overline{p},q)\in \mathbb{Z}^d\times \mathbb{N}$.
The number of solutions is given by 
$$
N_T(\overline{x})=|\{(\overline{p},q)\in \mathbb{Z}^d\times \mathbb{N}:\, 1\le q< T, \qen \hbox{ \eqref{eq:diop2} holds}  \}|.
$$
One can show, using Schmidt's arguments in \cite{sch1}, that for almost every $\overline{x}\in [0,1]^d$,
\begin{equation}
\label{eq:N_T}
N_T(\overline{x})= 2^dc_1\cdots c_d\, \log T+O_{x,\eps}((\log T)^{1/2+\eps})\quad \hbox{ for all $\eps>0$.}
\end{equation}
We prove here that the error term in \eqref{eq:N_T} satisfies the Central Limit Theorem. In the special case when all weights $(w_i)$ are equal, this result was established in \cite{dfv}.

\begin{theorem}[CLT for Diophantine approximation]
	\label{th:dioph}
	For $d\ge 2$, there exists explicit $\sigma>0$ such that for every $u\in \mathbb{R}$,
	$$
	\hbox{\rm Leb}\left(\left\{\overline{x}\in [0,1]^d:\, \frac{N_T(\overline{x}) - 2^dc_1\cdots c_d\, \log T }{(\log T)^{1/2}}<u \right\} \right)
	\longrightarrow
	\frac{1}{\sqrt{2\pi\sigma}} \int_{-\infty}^u e^{-t^2/(2\sigma)}\, dt
	$$
	as $T\to\infty$.
\end{theorem}

In this note we outline the proofs of Theorem \ref{th:lattice1}--\ref{th:dioph}.
Details will be published elsewhere. 

\section{Ingredients in the proofs}

\subsection{The space of lattices and Siegel transforms}
\label{subsec:siegel}
We denote by $\mathcal{L}_d$ the space of lattices in $\mathbb{R}^d$ with covolume one.
We recall that $\mathcal{L}_d$ can be realised as a homogeneous space
$\mathcal{L}_d\simeq \hbox{SL}_d(\mathbb{R})/\hbox{SL}_d(\mathbb{Z})$,
so that it is equipped with the unique invariant probability measure $\lambda_d$.
Given a bounded Borel measurable function $f:\mathbb{R}^d\to\mathbb{R}$ with compact support, we define its 
\emph{Siegel transform} $\widehat f:\mathcal{L}_d\to\mathbb{R}$ by
$$
\widehat f(\Lambda)=\sum_{v\in\Lambda\backslash\{0\}} f(v)\quad \hbox{ for $\Lambda\in \mathcal{L}_d.$}
$$

The starting point of our approach is the observation that the counting functions in Theorems \ref{th:lattice1}--\ref{th:dioph} can be realized as certain averages of suitable Siegel transforms. This idea is simpler to explain in the setting of Theorem \ref{th:dioph},
so let us begin by focusing on this case. Let 
$$
a=\hbox{diag}(2^{w_1},\ldots, 2^{w_d},2^{-1})\in\hbox{SL}_{d+1}(\mathbb{R})
$$
 and let $\chi$ denote the characteristic function of the domain
$$
\big\{(\overline{x},y)\in\mathbb{R}^d\times\mathbb{R} \, : \, 1\leq y < 2 \qen \hbox{\eqref{eq:diop2} holds}\big\}.
$$
Then one can readily check that for every $T=2^N$ with $N\ge 1$,
\begin{equation}
\label{eq:approx1}
N_{T}(\overline{x})=\sum_{n=0}^{N-1}\widehat \chi(a^n \Lambda_{\overline{x}}),
\end{equation}
where the lattice $\Lambda_{\overline{x}}$ is defined in \eqref{eq:lll}.
This basic formula allows to study the distribution of $N_{T}(\overline{x})$
using dynamics on the space of lattices. More precisely, we shall use an approximation 
of the form \eqref{eq:approx1} with $\chi$ replaced by a smooth function $f_\eps$ that approximates $\chi$
well in the $L^1$- and $L^2$-senses. \\

The counting functions in Theorem \ref{th:lattice1} and \ref{th:spiralling} can also be approximated along similar lines,
but the formulas are more complicated; in particular, one-parameter subgroups of diagonal matrices are no longer enough 
to achieve an approximation of $N_T$ as in \eqref{eq:approx1}. Let $A_d$ denote the subgroup of diagonal matrices in 
$\SL_d(\bR)$, and set $\theta_r=\hbox{diag}(1,\ldots,1,e^r)$. We shall show that for suitably chosen smooth compactly 
supported functions ${f}_{\eps,T}$ on $\mathbb{R}^d$ and finite subsets $B(r,T)$ of $A_d$,
\begin{equation}
\label{eq:approx2}
|\Lambda\cap \Omega_T|\approx \int_{\log a}^{\log b} \left(\sum_{a\in B(r,T)} \widehat{f}_{\eps,T}(\theta_r a\Lambda) \right)\, dr,
\end{equation}
in the $L^1$- and $L^2$-norms for $(\cL_d,\lambda_d)$. Our arguments from now on depend crucially on the fact (which will be explained in more detail below) that for a smooth compactly supported function $\phi$ on $\cL_d$, the collections of functions 
\[
\big\{ \phi(a \Lambda) \, : \, a \in A_d \big\}
\]
are "weakly independent".

\subsection{The method of cumulants}

There exists today a plethora of different techniques to establish convergence to the Gaussian distribution.
One of the first such techniques - if not the first - is nowadays often referred to as the \emph{Method of Moments}, 
and was used by Chebyshev to prove the classical Central Limit Theorem. We refer to \cite{bel} for a modern exposition 
of this technique. An essentially equivalent technique, but better tailored for problems pertaining to Gaussian distributions,
was later developed by Fr\'echet and Shohat, and goes under the name "Method of Cumulants". Let us briefly survey this
method. Given bounded random variables $X_1,\ldots,X_r$, their \emph{joint cumulant} is defined by the curious expression
$$
\cum^{(r)}(X_1,\ldots,X_r) = \sum_{\cP} (-1)^{|\cP|-1} (|P|-1)!\prod_{I \in \cP} 
\mathbb{E}\Big(\prod_{i\in I} X_i\Big),
$$
where the sum is taken over all partitions $\cP$ of the set $\{1,\ldots,r\}$.
We also set
$$
\cum^{(r)}(X) = \cum^{(r)}(X,\ldots,X)
$$
for a single bounded random variable $X$.
The cumulants have many useful combinatorial properties (see \cite{speed}).
For instance, if there exists a non-trivial partition
$\{1,\ldots,r\}=I\sqcup J$ such that the collections $\{X_i:\, i\in I\}$ and $\{X_j:\, j\in J\}$
are independent of each other, then
\begin{equation}
\label{zero}
\cum^{(r)}(X_1,\ldots,X_r)=0.
\end{equation}
Furthermore, a bounded random variable $X$ with mean zero is normally distributed if and only if $\cum^{(r)}(X)=0$ for $r\ge 3$.
In what follows, we shall use the following useful criterion due to Fr\'echet and Shohat \cite{FS} to establish our Central Limit Theorems.

\begin{proposition}[Method of Cumulants]
	\label{th:CLT}
	Let $(Z_T)$ be a collection of real-valued bounded random variables with mean zero satisfying
	\begin{align}
	&\sigma^2 := \lim_{T\to\infty} \hbox{\rm Var}(Z_T) < \infty,\label{eq:cond1}
	\end{align}
	and
	\begin{align}
	&\lim_{T\to\infty} \cum^{(r)}(Z_T) = 0, \quad \textrm{for all $r \geq 3$}.\label{eq:cond2}
	\end{align}
	Then for every $u\in\mathbb{R}$, 
	$$
	\hbox{\rm Prob}(Z_T<u)	\longrightarrow
	\frac{1}{\sqrt{2\pi\sigma}} \int_{-\infty}^u e^{-t^2/(2\sigma)}\, dt\quad \hbox{	as $T\to\infty$.}
	$$
\end{proposition}

In our recent work \cite{BG}, we used the Method of Cumulants to establish 
a general Central Limit Theorem for group actions which are exponentially 
mixing of all orders.  Here we essentially follow the approach developed in \cite{BG},
but substantial modifications will have to be made in order to handle more general averaging schemes,
as well as unbounded test functions.

\subsection{Estimates on the higher-order correlations}\label{sec:corr}
Let us now prepare the asymptotic formulas for higher-order correlations that will be used to 
estimate the cumulants and the variance. These formulas will be formulated in terms of Sobolev norms
$S_k$, $k\ge 1$, defined for smooth compactly supported functions on the space $\cL_d$ (see \cite{BEG}).
In the proofs of Theorems \ref{th:lattice1} and \ref{th:spiralling}, 
we will use estimates on correlations for the action on $\cL_d$ of the group of diagonal matrices 
$A_d\subset\hbox{SL}_d(\mathbb{R})$. If we fix a invariant metric $\rho$ on $A_d \cong \bR^{d-1}$, then the 
following result is a special case of \cite[Th.~1.1]{BEG}.

\begin{theorem}[Exponential multiple mixing of all orders]
	\label{th:mixing}
	For every $r\ge 2$, there exists an integer $k_r$ such that for all $k\ge k_r$, there is $\delta_{r,k}>0$ with the property that for 
	all $\phi_1,\ldots,\phi_r\in C_c^\infty(\mathcal{L}_d)$ and ${a}_1,\ldots,{a}_r\in A_d$, 
	\begin{align*}
	\int_{\mathcal{L}_d} \phi_1(a_1\Lambda)\cdots \phi_r(a_r\Lambda)\, d\lambda_d(\Lambda)=&\left(\int_{\mathcal{L}_d} \phi_1\, d\lambda_d\right)\cdots
	\left(\int_{\mathcal{L}_d} \phi_r\, d\lambda_d\right)\\
	&+ O_{r,k}\left( e^{-\delta_{r,k} D({a}_1,\ldots,{a}_r)}\, S_k(\phi_1)\cdots S_k(\phi_r)\right),
	\end{align*}
	where
	$D({a}_1,\ldots,{a}_r)=\min\{\rho(a_i,a_j):\, i\ne j\}.$
\end{theorem}

In order to study weighted Diophantine approximation \eqref{eq:diop2}, we need to analyze the distribution of orbits for 
the one-parameter semigroup
$$
a_w(t)=\hbox{diag}(e^{w_1 t},\ldots, e^{w_d t}, e^{-t} ),\quad t > 0,
$$
for lattices contained in the subset
$$
\mathcal{Y}_d=\left\{\Lambda_{\overline{x}}: \overline{x}\in [0,1]^d\right\}
$$
of $\cL_{d+1}$. We denote by $\sigma_d$ the measure $\cY_d$ induced by the Lebesgue measure on 
$[0,1]^d$. We establish the following asymptotic formula for the  higher-order correlations of the measures 
$a_w(t)_*\sigma_d$, generalizing the work of Kleinbock and Margulis \cite{km3}.

\begin{theorem}
\label{th:restricted}	
	For every $r\ge 2$ and $k\ge k_r$, there exists $\delta'_{r,k}>0$ such that for every $\phi_1,\ldots,\phi_r\in C_c^\infty(\mathcal{L}_{d+1})$ and ${t}_1,\ldots,{t}_r>0$, 
	\begin{align*}
	\int_{\mathcal{Y}_d} \phi_1(a_w({t}_1)y)\cdots \phi_r(a_{w}({t}_r)y)\, d\sigma_d(y)=&\left(\int_{\mathcal{L}_{d+1}} \phi_1\, d\lambda_{d+1}\right)\cdots
	\left(\int_{\mathcal{L}_{d+1}} \phi_r\, d\lambda_{d+1}\right)\\
	&+ O_{r,k}\left( e^{-\delta'_{r,k} D'({t}_1,\ldots,{t}_r)}\, S_k(\phi_1)\cdots S_k(\phi_r)\right),
	\end{align*}
	where $D'({t}_1,\ldots,{t}_r)=\min\{t_i, |{t}_i-{t}_j|:\, i\ne j\}.$
\end{theorem}

We note that $\cY_d$ is an unstable manifold for the one-parameter semigroup 
$$
g(t)=\hbox{diag}(e^{t/d},\ldots,e^{t/d},e^{-t}), \quad t>0,
$$
but not for semigroup $a_w(t)$, unless the weights $w_i$ are all equal.
The quantitative equidistribution for the translates $a_w(t)\cY_d$ has been established 
by Kleinbock and Margulis in \cite{km3}. We refine their argument to deal with higher-order correlations. The proof of Theorem \ref{th:restricted} was inspired by \cite{km3}.
It goes by induction on $r$
and uses quantitative equidistribution of the measures $a_{w}(t)_*\sigma_d$ combined with non-divergence estimates for the unipotent flows. \\

From Theorem \ref{th:restricted}, we deduce the following non-divergence estimate: 

\begin{corollary}[Non-divergence]
	\label{cor:non_div}
	Let $f$ be a continuous compactly supported function on $\mathbb{R}^d$.
	Then there exists $c>0$ such that for every $L\ge 1$ and $n\ge c\, \log L$,
	$$
	\sigma_d\left(\{y\in \cY_d:\, \widehat f(a^ny)>L \}\right)\ll_f L^{-d-1}.
	$$
\end{corollary}

This corollary will be used to construct bounded approximations for Siegel transforms.

\subsection{Bounded  approximations}\label{sec:bounded}

It might now be tempting to try to apply Proposition \ref{th:CLT},
combined with Theorems \ref{th:mixing} and \ref{th:restricted},
to the approximations \eqref{eq:approx1} and \eqref{eq:approx2} directly.
However,
we stress that the Siegel transform $\widehat f$ of a smooth compactly supported function
$f$ on $\mathbb{R}^d$ gives an \emph{unbounded} function on the space of lattices $\cL_d$. 
Moreover, the Sobolev norms $S_k(f)$ in \S\ref{sec:corr} are infinite.
In order to deal with these issues, we shall use that $\widehat f \in L^p(\lambda_d)$ for $p<d$ and 
show that one can approximate $\widehat f$ by a family of functions $\phi_{L}\in C_c^\infty(\mathcal{L}_d)$
satisfying
\begin{align}
\|\phi_L\|_\infty&=O(L) \qand S_k(\phi_L)=O_k(L^{d+1}),\label{eq:est1}
\end{align}
and
\begin{align}
\|\widehat{f} -\phi_L\|_1&=O_q (L^{-q})\;\hbox{ for all $q<d-1$},\label{eq:est2} 
\end{align}
and
\begin{align}
\|\widehat{f} -\phi_L\|_2&=O_q (L^{-q})\;\hbox{ for all $q<(d-2)/2$}. \label{eq:est3}
\end{align}
This observation will allow us to exploit the estimates from Subsection \ref{sec:corr} to analyze
the variance and the cumulants of higher orders. In the setting of Theorems \ref{th:lattice1} and \ref{th:spiralling}, 
we shall use the approximation \eqref{eq:approx2} and consider
\begin{align*}
Z_T(\Lambda) =
\int_{\log a}^{\log b} \left(\sum_{a\in B(r,T)} \widehat{f}_{\eps,T}(\theta_r a\Lambda) \right)\, dr,
\end{align*}
and
\begin{align*}
Z^*_T(\Lambda) &=
\int_{\log a}^{\log b} \left(\sum_{a\in B(r,T)} \phi^L_{\eps,T}(\theta_r a\Lambda) \right)\, dr,
\end{align*}
where $\phi^L_{\eps,T}$ is the bounded approximation for $\hat{f}_{\eps,T}$. Because of \eqref{eq:est2}, the parameter 
$L=L(T)$ can be chosen so that 
$$
\|Z_T-Z^*_T\|_{L^1(\lambda_d)}\to 0\quad\hbox{as $T\to\infty$.}
$$
After this choice has been made, it suffices to analyze convergence in distribution of $Z_T^*$. \\

In the proof of Theorem \ref{th:dioph}, we consider
$$
Z_T(y)=\frac{1}{\sqrt{N}} \sum_{n=0}^{N-1} \left(\widehat{\chi}(a^n y) -\int_{\mathcal{Y}_d}\widehat{\chi}(a^n y)\, d\sigma_d(y)\right)
$$
and its approximation 
$$
Z^*_T(y)=\frac{1}{\sqrt{N}} \sum_{n=0}^{N-1} \left(\phi_{\eps}^L(a^n y) -\int_{\mathcal{Y}_d}\phi_{\eps}^L(a^n y)\, d\sigma_d(y)\right),
$$
where $\phi_{\eps}^L$ denotes the bounded approximation to $\widehat f_\eps$ as above, and 
$f_\eps$ is a smooth approximation to the characteristic function $\chi$ (recall the notation from Subsection \ref{subsec:siegel}). 
Here the parameters $\eps=\eps(T)$ and $L=L(T)$ can be chosen so that
\begin{equation}
\label{eq:Z}
\|Z_T-Z^*_T\|_{L^1(\sigma_d)}\to 0\quad\hbox{as $T\to\infty$.}
\end{equation}
To arrange \eqref{eq:Z}, we use the non-divergence estimate established in Corollary \ref{cor:non_div} and the following uniform bound
\begin{equation}
\label{eq:uni}
\sup_{n\ge 1} \left\|\hat f_\eps \circ a^n\right\|_{L^2(\sigma_d)}<\infty.
\end{equation}
In order to prove \eqref{eq:uni}, we interpret the $L^2$-norm arithmetically and 
reduce this estimate to a problem of counting solutions of certain Diophantine equations. \\

Ultimately, we shall show that $Z^*_T$ converges to the Normal Law using Proposition \ref{th:CLT}.
Our main tool is the estimates on higher-order correlations from \S\ref{sec:corr}.
We note the bounds in our computations will depend on the parameters $\eps,T,L$, and thus
the explicit forms of error terms in Theorems \ref{th:mixing} and \ref{th:restricted}
are essential for this purpose.

\subsection{Well-separated tuples and estimating the cumulants}
By linearity, the estimates on cumulants arising in the proofs of Theorems \ref{th:lattice1} and
\ref{th:spiralling}	 reduce to the following basic problem, which is discussed in more detail in 
our paper \cite{BG}. Given $\phi_1,\ldots,\phi_r\in C_c^\infty(\mathcal{L}_d)$ and $(a_1,\ldots,{a}_r)\in A_d^r$, we 
wish to estimate averages of cumulants of the form
\begin{equation}
\label{eq:c1}
\cum^{(r)}_{\lambda_d}\left(\phi_1\circ a_1,\ldots, \phi_r\circ a_r\right)
= \sum_{\cP} (-1)^{|\cP|-1} (|P|-1)!\prod_{I \in \cP} \left(
\int_{\cL_d} \Big(\prod_{i\in I} \phi_i\circ a_i\Big)\,d\lambda_d\right),
\end{equation}
as $(a_1,\ldots,a_r)$ varies over certain subsets of $A_d^r$. \\

The idea is to decompose $A_d^r$ into finitely many regions where the cumulants can be  estimated separately. These regions are defined as follows. Recall that $\rho$ is a fixed invariant metric on $A_d$. For $I, J \subset [r]$ and $\overline{a} = (a_1,\ldots,a_r) \in A_d^r$,
we set
$$
\rho^{I}(\overline{a}) = \max\big\{ \rho(a_i,a_j) \, : \, i,j \in I \big\}\quad\hbox{and}
\quad
\rho_{I,J}(\overline{a}) = \min\big\{ \rho(a_i,a_j) \, : \, i \in I, \enskip j \in J \big\}.
$$
If $\cQ$ is a partition of $\{1,\ldots,r\}$, we define
$$
\rho^{\cQ}(\overline{a}) = \max\big\{ \rho^{I}(\overline{a}) \, : \, I \in \cQ \big\}
\quad\hbox{and}\quad
\rho_{\cQ}(\overline{a}) = \min\big\{ \rho_{I,J}(\overline{a}) \, : \, I \neq J, \enskip I, J \in \cQ \big\}.
$$
For $0 \leq \alpha < \beta$, we define
$$
\Delta_{\cQ}(\alpha,\beta) 
= 
\big\{ 
\overline{a} \in A_d^r \, : \, 
\rho^{\cQ}(\overline{a}) \leq \alpha, 
\qen
\rho_{\cQ}(\overline{a}) > \beta \big\}
$$
and
$$
\Delta(\beta)= 
\big\{ 
\overline{a} \in A_d^r \, : \, 
\rho(a_i,a_j) \leq \beta \hbox{ for all $i,j$}\big\}.
$$
We shall think of the tuples in $\Delta_{\cQ}(\alpha,\beta)$, for some partition $\cQ$ with $|Q|\ge 2$, as being
``well-separated'', 
while we think of the tuples in $\Delta(\beta)$ as being ``clustered''.  \\

We estimate the cumulants on $\Delta_{\cQ}(\alpha,\beta)$ using the exponential multiple mixing property established in Theorem \ref{th:mixing}. 

\begin{lemma}[Proposition 6.1, \cite{BG}]
	\label{l:cum}
	For all $r\ge 3$ and $k>k_r+r$,
	there exist $c_{r,k},  \delta_{r,k}>0$ such that 
	for any partition $Q$ of $\{1,\ldots,r\}$ with $|Q|\ge 2$ and $s>0$,
	 we have
	$$
	|\cum^{(r)}_{\lambda_d}\left(\phi_1\circ a_1,\ldots, \phi_r\circ a_r\right)|\ll_{r,k}\, e^{-\delta_{r,k} s}\, S_k(\phi_1)\cdots S_k(\phi_r)
	$$
	when $(a_1,\ldots,a_r)\in \Delta_{\cQ}(s,c_{r,k}\,s)$.
\end{lemma}

To prove this lemma, we introduce a cumulants ``conditioned'' on a given partition $Q$.
For partitions $\cP$ and $\cQ$, we set $\cP\wedge \cQ=\{P\cap Q:\, P\in \cP,Q\in \cQ\}$.
We define
\begin{equation}
\label{eq:c2}
\cum^{(r)}_{\lambda_d,\cQ}\left(\phi_1\circ a_1,\ldots, \phi_r\circ a_r\right)
= \sum_{\cP} (-1)^{|\cP|-1} (|P|-1)!\prod_{J \in \cP\wedge \cQ} 
\left(\int_{\cL_d} \Big(\prod_{i\in J} \phi_i\circ a_i\Big)\,d\lambda_d\right).
\end{equation}
Comparing \eqref{eq:c1} and \eqref{eq:c2}, we realize that they are approximately equal
for tuples $(a_1,\ldots,a_r)\in \Delta_{\cQ}(\alpha,\beta)$ with suitably chosen $\alpha$ and $\beta$ because
$$
\int_{\cL_d} \Big(\prod_{i\in I} \phi_i\circ a_i\Big)\,d\lambda_d
\approx
\prod_{K \in \cQ} 
\int_{\cL_d} \Big(\prod_{i\in I\cap K} \phi_i\circ a_i\Big)\,d\lambda_d, \quad \textrm{for all $I, K \subset [r]$},
$$
according to Theorem \ref{th:mixing}.
The second step in the proof of Lemma \ref{l:cum} utilizes the fact that
when $Q$ is a non-trivial partition of $[r]$, then
$$
\cum^{(r)}_{\lambda_d,\cQ}\left(\phi_1\circ a_1,\ldots, \phi_r\circ a_r\right)
=0,
$$
which is a combinatorial version of \eqref{zero} (see Proposition 8.1 in \cite{BG}).
This leads to the estimate in Lemma \ref{l:cum}.

\vspace{0.2cm}

In order to apply Lemma \ref{l:cum}, we decompose $A_d^r$ into regions
where the tuples $({a}_1,\ldots,{a}_r)$ are ``well-separated'' or ``clustered'' on certain scales. 
We show (cf. \cite[Prop. 6.2]{BG}) that for suitably chosen parameters
$$
0=\alpha_0<\beta_0<\alpha_1<\cdots<\beta_{r-1}<\alpha_{r},
$$
we have a decomposition
\begin{equation}
\label{eq:decomp}
	A_d^r = \Delta(\alpha_r) \cup \Big( \bigcup_{j=0}^{r-1} \bigcup_{|\cQ| \geq 2} \Delta_{\cQ}(\alpha_j,\beta_{j}) \Big),
\end{equation}
where the union is taken over the partitions $\cQ$ of $\{1,\ldots,r\}$ with $|Q|\ge 2$.
It turns out possible to choose the parameters $\alpha_j,\beta_j$ in such a way
that Lemma \ref{l:cum} can be applied to the averages of the cumulants over subsets of
$\Delta_{\cQ}(\alpha_j,\beta_{j})$ to conclude that they are negligible.
Now it remains to estimate
the average over a subset of $\Delta(\alpha_r)$. Since we can choose $\alpha_r$ quite small,
the latter average can be estimated by bounding the number of terms.

\vspace{0.2cm}

The above argument requires some modifications for the proof of Theorem \ref{th:dioph}
because we need to take into account the estimator $D'$ in Theorem \ref{th:restricted}.
It will be convenient to embed $A_d^r$ in $A_d^{r+1}$
by $\overline{a}\mapsto (e,\overline{a})$ and define subsets $\Delta_\cQ(\alpha,\beta)$
of $A_d^r$
with respect to this embedding
for partitions $\cQ$ of $\{0,1,\ldots,r\}$.
As before, we use the decomposition \eqref{eq:decomp}.
When the partition $\cQ$ is non-trivial and different from $\{\{0\},\{1,\ldots,r\}\}$,
we are able to modify the proof of Lemma \ref{l:cum} using Theorem \ref{th:restricted}
and estimate the cumulants $\cum^{(r)}_{\sigma_d}\left(\phi_1\circ a^{n_1},\ldots, \phi_r\circ a^{n_r}\right)$ when $(a^{n_1},\ldots,a^{n_r})\in \Delta_{\cQ}(s,c_{r,k}\,s)$.
When $Q=\{\{0\},\{1,\ldots,r\}\}$,
we observe that Theorem \ref{th:restricted} implies that 
$$
\cum^{(r)}_{\sigma_d}\left(\phi_1\circ a^{n_1},\ldots, \phi_r\circ a^{n_r}\right)
\approx
\cum^{(r)}_{\lambda_d}\left(\phi_1\circ a^{n_1},\ldots, \phi_r\circ a^{n_r}\right)
$$
when $(a^{n_1},\ldots,a^{n_r})\in \Delta_{\cQ}(s,c_{r,k}\,s)$,
and the latter cumulant has already been estimated. \\

Finally, we have to deal with the average over $(a^{n_1},\ldots,a^{n_r})\in\Delta(\alpha_r)$.
For this purpose, we modify the function $Z^*_T$ in such a way that its convergence in distribution is not affected.
We set 
$$
Z^{**}_T(y)=\frac{1}{\sqrt{N}} \sum_{n=M}^{N-1} \left(\phi_{\eps}^L(a^n y) -\int_{\mathcal{Y}_d}\phi_{\eps}^L(a^n y)\, d\sigma_d(y)\right),
$$
where the parameter $M=M(N)\to \infty$  is chosen so that
$$
\|Z^{*}_T-Z^{**}_T\|_{L^1(\sigma_d)}\to 0\quad\hbox{as $T\to\infty$.}
$$
In particular $Z^*_T$ and $Z_T^{**}$ have the same distributional limits, and thus it is suffices to establish 
convergence in distribution of $Z^{**}_T$. Choosing $M = M(N)$ appropriately, we can further arrange so 
that averages over subsets of $\Delta(\alpha_r)$ in the cumulant calculations  $\cum_{\sigma_d}^{(r)}(Z^{**}_T)$ 
tend to zero.

\subsection{Estimating the variance}
In order to estimate the variance in the setting of Theorems \ref{th:lattice1} and \ref{th:spiralling},
we need to consider sums of the form
\begin{equation}
\label{eq:var1}
\frac{1}{|B_T|}\sum_{a,b\in B_T}\int_{\cL_d} \psi (a\Lambda)\psi (b\Lambda)\, d\lambda_d(\Lambda),
\end{equation}
where $B_T$'s are finite subsets of a lattice $\Delta_d\subset A_d$ and $\psi=\phi-\int_{\cL_d}\phi\, d\lambda_d$
with a smooth compactly supported function $\phi$ on $\cL_d$ (compare with \eqref{eq:approx1}).
Using invariance of the measure $\lambda_d$, we rewrite this expression as
$$
\frac{1}{|B_T|}\sum_{a\in B_T, c\in a^{-1}B_T}\int_{\cL_d} \psi (a\Lambda)\psi (ac\Lambda)\, d\lambda_d(\lambda)
=\sum_{c\in B_T^{-1}B_T} \frac{|B_T\cap B_Tc^{-1}|}{|B_T|}
\int_{\cL_d} \psi (\psi\circ c)\, d\lambda_d.
$$
Then using Theorem \ref{th:mixing} with $r=2$, we deduce that \eqref{eq:var1}
converges to 
\begin{equation}
\label{eq:var1_0}
\sum_{c\in \Delta_d} 
\int_{\cL_d} \psi (\psi\circ c)\, d\lambda_d.
\end{equation}
It should be noted that this argument have to be applied to the family of functions
$\phi_{\eps,T}^L\circ \theta_r$, introduced in \S\ref{sec:bounded},
with suitably chosen parameters.
The explicit form of the error term in Theorem \ref{th:mixing}
still allows to justify convergence of \eqref{eq:var1}.

\vspace{0.2cm}

The computation of variance in Theorem \ref{th:dioph} reduces to analysing the expressions
\begin{align}
\frac{1}{N}\sum_{n,m=0}^{N-1}\left(\int_{\cY_d} \psi(n,y)\psi(m,y)d\sigma_d(y)\right)\label{eq:var2}
=&
\frac{1}{N}\sum_{n=0}^{N-1}\left(\int_{\cY_d} \psi(n,y)^2\,d\sigma_d(y)\right)\\
&+\frac{2}{N}\sum_{0\le n<m\le N-1}\left(\int_{\cY_d} \psi(n,y)\psi(m,y)d\sigma_d(y)\right),\nonumber
\end{align}
where $\psi(n,y)=\phi(a^n y)-\int_{\cY_d} (\phi\circ a^n) \, d\sigma_d$ for smooth compactly supported functions $\phi$ on $\cL_{d+1}$.
We shall show that \eqref{eq:var2} converges to 
\begin{equation}
\label{eq:var2_0}
\sum_{k\in \mathbb{Z}}\left(\int_{\cL_{d+1}} \phi(\phi\circ a^k) \, d\lambda_{d+1}-\left(\int_{\cL_{d+1}} \phi \,d\lambda_{d+1}\right)^2\right)
\end{equation}
as $N\to\infty$.
First, we observe that by Theorem \ref{th:restricted} the first term 
in \eqref{eq:var2} converges
to 
$$
\int_{\cL_{d+1}} \phi^2 \, d\lambda_{d+1}-\left(\int_{\cL_{d+1}} \phi \,d\lambda_{d+1}\right)^2
$$
as $N\to\infty$.
To estimate the second term in \eqref{eq:var2}, we rewrite it as
\begin{equation}
\label{eq:var2_1}
\sum_{k=1}^{N-1} \left(\frac{2}{N} \sum_{n=0}^{N-1-k} \left(\int_{\cY_d} \psi(n,y)\psi(n+k,y)d\sigma_d(y)\right) \right).
\end{equation}
It follows from Theorem \ref{th:restricted} that for fixed $k$,
$$
\int_{\cY_d} \psi(n,y)\psi(n+k,y)d\sigma_d(y)\longrightarrow
\int_{\cL_{d+1}} \phi(\phi\circ a^k) \, d\lambda_{d+1}-\left(\int_{\cL_{d+1}} \phi \,d\lambda_{d+1}\right)^2
$$
as $n\to \infty$. A more tedious analysis, which utilizes the explicit quantitative bounds
from Theorem \ref{th:restricted} with $r=1,2$, allows to conclude that \eqref{eq:var2_1}
converges to 
$$
2\sum_{k\ge 1}\left(\int_{\cL_{d+1}} \phi(\phi\circ a^k) \, d\lambda_{d+1}-\left(\int_{\cL_{d+1}} \phi \,d\lambda_{d+1}\right)^2\right).
$$ 
This leads to the formula \eqref{eq:var2_0}.
More precisely, this argument will be applied to the family of functions $\phi_\eps^L$, 
introduced in \S\ref{sec:bounded},
but the explicit form of the error term in Theorem \ref{th:restricted}
allows to handle this.

Finally, we note that the expressions \eqref{eq:var1_0} and \eqref{eq:var2_0} can be computed
explicitly in our setting using Rogers' formula \cite{rog}. In particular, we conclude that
the obtained variances are positive.

\end{document}